\documentstyle[12pt]{article}
\newlength{\abstractwidth}
\renewcommand{\thefootnote}{\fnsymbol{footnote}}
\renewcommand{\thanks}[1]{\footnote{#1}} % Use this for footnotes
\newcommand{\starttext}{ \setcounter{footnote}{0}
\renewcommand{\thefootnote}{\arabic{footnote}}}

\newcommand{\be}{\begin{equation}}
\newcommand{\bea}{\begin{eqnarray}}
\newcommand{\eea}{\end{eqnarray}} \newcommand{\ee}{\end{equation}}
 
 \def\ba{\begin{eqnarray}}
\def\ea{\end{eqnarray}}

\def\o{\omega}

\def\tr{{\rm tr}}
\def\det{{\rm det}}

\def\log{\,{\rm log}\,}
\def\exp{\,{\rm exp}\,}

\def\o{\omega}

\def\o{\omega}

\def\ge{\geq}
\def\le{\leq}

\def\p{\partial}

\def\Ent{{\rm Ent}}

\def\[{{\bf [}}
\def\]{{\bf ]}}

\def\ddbar{i\p\bar\p}

\def\mathbb{\bf}

\def\eqref{\ref}

\begin{document}
\starttext \baselineskip=18pt \setcounter{footnote}{0}
\newtheorem{theorem}{Theorem}
\newtheorem{lemma}{Lemma}
\newtheorem{corollary}{Corollary}
\newtheorem{definition}{Definition}
\newtheorem{conjecture}{Conjecture}
\newtheorem{proposition}{Proposition}

\begin{center}
{\Large \bf ON $L^\infty$ ESTIMATES FOR COMPLEX MONGE-AMP\`ERE EQUATIONS
\footnote{Work supported in part by the National Science Foundation under grant DMS-1855947.}}

\medskip
\centerline{Bin Guo, Duong H. Phong, and Freid Tong}

\medskip

\begin{abstract}

{\footnotesize A PDE proof is provided for the sharp $L^\infty$ estimates for the complex Monge-Amp\`ere equation which had required pluripotential theory before. The proof covers both cases of fixed background as well as degenerating background metrics. It extends to more general fully non-linear equations satisfying a structural condition, and it also gives estimates of Trudinger type.}

\end{abstract}

\end{center}

\baselineskip=15pt
\setcounter{equation}{0}
\setcounter{footnote}{0}

\section{Introduction}
\setcounter{equation}{0}

The main goal of this paper is to answer a long-standing question in the theory of complex Monge-Amp\`ere equations and its applications to complex geometry, namely whether  sharp $L^\infty$ estimates can be established by PDE methods and without pluripotential theory. We shall see that the answer is affirmative and, as may have been anticipated, the PDE proof also gives new estimates as well as extensions to many other fully non-linear equations.

\smallskip
$L^\infty$ estimates have a uniquely storied history in the theory of complex Monge-Amp\`ere equations. Early on, they were recognized as the defining difficulty in the problem of finding K\"ahler-Einstein metrics. Yau's introduction in 1978 of Moser iteration was the key step in his solution of the Calabi conjecture, and it ushered in a new era for complex Monge-Amp\`ere equations \cite{Y}. Yau's Moser iteration method works for equations whose right hand sides are in $L^q$ for $q>n$, where $n$ is the complex dimension of the underlying space. The next major advance was the 1998 result of Kolodziej \cite{K}, which established $L^\infty$ estimates for the solution when the right hand side is in $L^q$, for $q>1$. This improvement in the range of $q$ is no mere technicality, and it has a profound geometric significance: that $L^\infty$ estimates fail for $q=1$ is indicative of the necessity of stability conditions in the K\"ahler-Ricci flow \cite{PSS}, while the cases $1<q\leq n$ are needed in a wide range of problems, including singular K\"ahler-Einstein metrics on manifolds of general type \cite{EGZ, DP}, the analytic minimal model program \cite{ST}, and degenerating Calabi-Yau metrics \cite{T}. Kolodziej's method of proof relied heavily on the pluripotential theory developed in the late 1970's by Bedford and Taylor \cite{BT1,BT2}. Actually, for many applications to geometry, an extension of Kolodziej's results to the more general case of degenerating background metrics is necessary. Such an extension was developed in 2007 independently by Demailly and Pali \cite{DP} and Eyssidieux, Guedj, and Zeriahi \cite{EGZ}, and pluripotential theory continued to be essential.

\smallskip
An immediate question for the theory of geometric partial differential equations is whether the above $L^\infty$ estimates can be derived by PDE methods, instead of pluripotential theory which is specific to Monge-Amp\`ere equations. This question gained considerable attention over the years, as more and more fully non-linear equations without corresponding pluripotential theory emerged in complex differential geometry.  A completely different proof for right hand sides in $L^q$ with $q> 2$ was found in 2011 by Blocki \cite{B11}, using the Alexandrov-Bakelman-Pucci (ABP) maximum principle.\footnote{The possibility of applying ABP maximum principle to the complex Monge-Amp\`ere equation were suggested a while ago by S.Y. Cheng and S.T. Yau.} His methods turn out to be remarkably powerful, and have since been applied successfully to many problems, including subsolutions \cite{Sze, PT}, equations with gradient terms \cite{TW}, and the constant scalar curvature problem \cite{CC}. Even so, extensions to degenerating backgrounds as well as the full range $q>1$ remained out of reach. In a different direction, a PDE proof of $L^\infty$ estimates for the complex Monge-Amp\`ere equation was obtained by J.X. Wang, X.J. Wang, and B. Zhou \cite{WWZ} for domains in ${\bf C}^n$. However, their methods do not appear adaptable to the compact manifold case, even in the simplest case when the background K\"ahler metric is fixed. Thus a fully effective approach to $L^\infty$ estimates remained an open question.

\smallskip
The PDE approach to $L^\infty$ estimates which we present in this paper combines the methods of Wang, Wang, and Zhou \cite{WWZ} with a fundamental new idea of Chen and Cheng \cite{CC} in their recent work on constant scalar curvature K\"ahler metrics, namely to compare the given equation with an auxiliary complex Monge-Amp\`ere equation. A key novelty in our paper resides in the choice of auxiliary complex Monge-Amp\`ere equation, as well as of the test function $\Phi$ for the comparison. We now formulate our main results.

\medskip
Let $X$ be a compact K\"ahler manifold without boundary of dimension $n$, and $\o_X$ its K\"ahler form. If $\varphi$ is a real smooth function on $X$, we let $\o_{\varphi}=\o_X+i\p\bar\p\varphi$, and let $h_{\varphi}$ be the corresponding endomorphism relative to the metric $\o_X$. Explicitly, if we write $\o_X=ig_{\bar kj}dz^j\wedge d\bar z^k$ in local holomorphic coordinates, then
$(h_{\varphi})^j{}_k=g^{j\bar m}(\o_{\varphi})_{\bar m k}$. Let $\lambda[h_{\varphi}]$ be the vector of eigenvalues of $h_{\varphi}$, and consider the non-linear partial differential equation
\bea
\label{eqn:main1}
f(\lambda[h_{\varphi}])=\,e^{F},
\quad
{\rm sup}_X\varphi=0,
\quad \lambda[h_{\varphi}]\in\Gamma,
\eea
for a given function $f(\lambda)$ and real function $F$ normalized such that $\int_X e^{nF}\omega_X^n = V = \int_X \omega_X^n$. Here the function $f(\lambda)$ is assumed to be invariant under permutations of the components of $\lambda$, and defined on an open cone $\Gamma\subset \{\lambda : \lambda_1 +\ldots + \lambda_n>0\}$ with vertex at the origin and containing the first octant
$\{\lambda: \lambda_1>0,\ldots, \lambda_n>0\}$. We assume throughout $f$ is elliptic in the sense that $\frac{\partial f(\lambda)}{\partial \lambda_j}>0$ for any $\lambda\in \Gamma$, and that
\bea\label{eqn:homo 1}
f(r\lambda)=r\,f(\lambda), \quad r>0,\ \lambda\in\Gamma,
\eea
which can be viewed as a normalization, as the same homogeneous equation can be expressed with many different functions $f(\lambda)$. We shall need the $L^1(\log L)^p$ norm of $e^{nF}$ with respect to the measure $\o_X^n$, which can be recognized as a generalized entropy $\Ent_p(F)$,
\bea
\label{eqn:entropy}
\Ent_p(F)={1\over V}\int_X e^{nF}|F|^p\o_X^n={1\over n^pV}\|e^{nF}\|_{L^1(\log L)^p}.
\eea

\begin{theorem}
\label{thm:main1}
Consider the equation (\ref{eqn:main1}), and assume that the function $f(\lambda)$ satisfies the following structural condition, namely that there exists a constant $\gamma>0$ so that
\bea\
\label{eqn:strf}
\det \big(\frac{\partial f(\lambda[h])}{\partial h_{ij}} \big) \ge \gamma ,\quad \mbox{for all }\lambda\in \Gamma. 
\eea
Fix $p >n$. Then for any solution $\varphi\in C^2(X)$, we have the estimate
\bea
{\rm sup}_X|\varphi|\leq C
\eea
where the constant $C$ depends only on $n,p,\gamma$, $\o_X$, and the entropy $\Ent_p(F)$.
\end{theorem}

We observe that many equations satisfy the structural condition (\ref{eqn:strf}). They include the Monge-Amp\`ere equation $f(\lambda)=(\prod_{j=1}^n\lambda_j)^{1\over n}$, more generally Hessian equations  $f(\lambda) = \sigma_k(\lambda)^{1/k}$, $k=1,\ldots, n$, where $\sigma_k(\lambda)$ is the $k$-th symmetric function, and quotient Hessian equations such as
$f(\lambda) = \big({\sigma_k\over\sigma_l} \big)^{1/(k-l)} + c \sigma_m^{1/m}$ for some $c>0$. Applying Theorem \ref{thm:main1} to the Monge-Amp\`ere equation, we obtain immediately Kolodziej's \cite{K} sharp $L^\infty$ bounds. Applying it to Hessian equations, we obtain the $L^\infty$ bounds of Dinew and Kolodziej \cite{DK}. In fact, our result is stronger, as the bound in \cite{DK} requires that $\exp(nF)$ be in $L^q$  for $q>1$, while we only need that $\exp(nF)$ be in $L^1(\log L)^p$ for $p>n$. Beyond these cases, the $L^\infty$ estimates in Theorem \ref{thm:main1} are new, and appear to be the first   obtained in any generality for fully non-linear equations from K\"ahler geometry.

\medskip
We discuss now estimates, particularly important for many geometric applications, where the background metric is allowed to degenerate. It is convenient to set up the equation as follows. Let $(X,\o_X)$ be a compact K\"ahler manifold of dimension $n$ as before, and let $\chi$ be a fixed closed and non-negative $(1,1)$-form. For $t\in (0,1]$, set
\bea
\o_t=\chi+t\o_X
\eea
and for each $\varphi\in C^2(X)$, let $\o_{t,\varphi}=\o_t+i\p\bar\p\varphi$, $h_{t,\varphi}$ be the corresponding endomorphism relative to the metric $\o_X$, and $\lambda[h_{t,\varphi}]$ be the vector of eigenvalues of $h_{t,\varphi}$. Consider the family of non-linear partial differential equations
\bea
\label{eqn:main}
f(\lambda[h_{t,\varphi_t}])=c_t\,e^{F_t},
\quad
{\rm sup}_X\varphi_t=0,
\quad \lambda[h_{t,\varphi_t}]\in\Gamma,
\quad t\in (0,1]
\eea
for given function $f(\lambda)$, real functions $F_t$, and positive constant coefficients $c_t$. Here the functions $F_t$ are normalized by $\int_X e^{nF_t} \omega_X^n = \int_X \omega_X^n$, so that the constants $c_t$ are determined. Denote by $V_t$ the volume of the metric $\o_t$, $V_t=\int_X\o_t^n$, and define the energy $E_t(\varphi)$ by
\bea
\label{eqn:E t}
E_t(\varphi)
= \frac{1}{V_t}\int_X (-\varphi) f(\lambda([h_{t,\varphi}])^n\omega_X^n.
\eea
For functions $\varphi_t$ solving the equation (\ref{eqn:main}), we obviously have
\bea
E_t(\varphi_t)= \frac{c^n_t}{V_t} \int_X (-\varphi_t) {\rm exp}(nF_t)\omega_X^n.
\eea

\begin{theorem}
\label{thm:main}
Consider the family of equations (\ref{eqn:main}), and assume that $f(\lambda)$ satisfies the structural condition (\ref{eqn:strf}). 
Let $\varphi_t$ be $C^2$ functions on $X$ satisfying (\ref{eqn:main}). Fix $p>n$. Then for any $t\in (0,1]$, we have
\bea
\sup_X|\varphi_t|\leq C
\eea
where $C$ is a constant depending only on $\o_X,\chi, p, n, \gamma$, and upper bounds for the following three quantities
\bea
\label{thm:key}
{c_t^n\over V_t},
\quad
E_t(\varphi_t),
\quad
\Ent_p(F_t).
\eea
\end{theorem}

All three quantities in (\ref{thm:key}) have attractive interpretations. We have already noted from (\ref{eqn:entropy}) that $\Ent_p(F_t)$ is a generalized entropy.
The quantity $c_t^n\over V_t$ can be viewed as a relative volume, and a substitute for a cohomological constraint when dealing with general equations $f(\lambda)$. The quantity $E_t(\varphi_t)$ is clearly an energy functional, as it reduces to the Dirichlet integral in the case of the Laplacian on surfaces.

\medskip

As a special case, Theorem \ref{thm:main} applied to the Monge-Amp\`ere equation gives back immediately the $L^\infty$ estimates of Eyssidieux, Guedj, Zeriahi \cite{EGZ}, and Demailly, Pali \cite{DP}, in the full generality of degenerating background metrics. The point here is that the relative volumes $c_t^n/V_t$ must be $1/V$ because of a cohomological constraint, and it follows easily from Jensen's inequality that the energies $E_t$ are uniformly bounded (see the fuller discussion in Theorem \ref{thm:MA} in \S 4 below). Thus Theorem \ref{thm:main} gives the pure PDE proof of these $L^\infty$ estimates that we sought. In fact, it is particularly simple as the ABP maximum principle is not even needed.

\smallskip

More generally, as long as the admissible cone $\Gamma$ is the one corresponding to PSH functions, we can obtain uniform bounds for the energies $E_t$. For example, this applies to the equations $f(\lambda) = \big( \frac{\sigma_k}{\sigma_l} \big)^{1/(k-l)} + c \; \sigma_n^{1/n}$. 

\medskip

For Hessian equations $f(\lambda)=\sigma_k(\lambda)^{1\over k}$, $1\leq k< n$, nothing was known in the case of degenerating background metrics, and Theorem \ref{thm:main} is now the only available result. As we shall see in Theorem \ref{thm:fully} in \S 5 below, if $\| e^{nF_t}\|_{L^p}$ is uniformly bounded for $p>{n\over k}$, we can bound the energy $E_t$ by a multiple of ${c_t^n V_t^{-1}}$. In particular, for big cohomology classes $[\o_t]$, the volumes $V_t$ do not tend to $0$, the coefficients $c_t^n$ are bounded by the entropy, and we obtain then uniform $L^\infty$ bounds. Such bounds are completely new for fully non-linear equations, and in particular Theorem \ref{thm:big} is new for Hessian equations.

\medskip

This paper is organized as follows. In \S 2 and \S 3, we give the proofs of Theorem \ref{thm:main1} and Theorem \ref{thm:main}. We actually begin with the proof of Theorem \ref{thm:main} in \S 2, as Theorem \ref{thm:main1} would follow from Theorem \ref{thm:main} upon control of the energy term $E_t$.  As in \cite{CC}, we use a comparison of the given equation to an auxiliary complex Monge-Amp\`ere equation. However, for Theorem \ref{thm:main}, it is very important to choose the auxiliary equation so as to avoid having to use the ABP maximum principle, as this maximum principle would be an impediment in the case of degenerating background metrics. The additional step of estimating the energy terms to get Theorem \ref{thm:main1} from Theorem \ref{thm:main} is provided by Theorem \ref{thm:energy}, which is the one requiring the ABP maximum principle. In \S 4 and \S 5, we show how our results apply to Monge-Amp\`ere and Hessian equations respectively, and how they improve on many recent results in the literature. Finally, 
we have discussed exclusively so far $L^\infty$ estimates. But not surprisingly, the same methods apply to $L^p$ and exponential estimates as well. We illustrate this in \S 6 by some applications to Trudinger exponential inequalities for fully non-linear equations, which are either new or independent proofs of known sharp estimates. We also observe that the assumption (\ref{eqn:homo 1}) that $f(\lambda)$ was homomegenous of degree $1$ was only for simplicity. The results of this paper still hold with this assumption replaced by the weaker Euler inequality with some fixed positive constant $\Lambda$,
\bea\nonumber
\sum_{j=1}^n{\p f\over\p \lambda_j}\lambda_j\leq \Lambda\,f(\lambda),
\quad\lambda\in \Gamma.
\eea
\noindent
Furthermore, our methods can be  adapted to the setting of families of K\"ahler manifolds $(X_j, \omega_j)$ of the same dimension, as long as the $\alpha$-invariant estimates hold uniformly.

\section{Proof of Theorem \ref{thm:main}}
\setcounter{equation}{0}
\renewcommand{\eqref}[1]{(\ref{#1})}

Let $\varphi_t$ solve the equation \eqref{eqn:main}. Fix $p>n$, and any upper bound ${\overline E}_t>0$ for $E_t(\varphi_t)$. We shall actually show that
$$
\sup_X |\varphi_t|\le C_0\Big\{ \frac{c_t^n}{V_t}  \big(\Ent_p + 1 +  {\rm exp}(C_1 {\overline E}_t)\big)\Big\}^{\frac{n}{p-n}} {\overline E}_t+ C(n,p),
$$
for constants $C_0, C_1$ depending only on $n,\omega_X,\chi, \gamma$, and $C(n,p)$ depending only on $n,p$.
Throughout the proof we will fix $t\in (0,1]$, but the constants will be independent of $t$, unless stated explicitly otherwise.

\medskip

For any $s>0$, we let $\Omega_s : = \{\varphi_t \le - s\}$ be the sub-level set of $\varphi_t$. 
\begin{lemma}
 \label{lemma trudinger} 
 There are constants $C= C(n, \omega_X,\chi,\gamma)>0$ and $\beta_0=\beta_0(n, \omega_X, \chi,\gamma)>0$ such that for any $s>0$
  \begin{equation}%\label{eqn:1.10}
\nonumber \int_{\Omega_s} \exp\Big\{  \beta_0 \big(\frac{- (\varphi_t + s)}{ A_{s}^{1/(n+1)}  } \big)^{\frac{n+1}{n}}   \Big\}\omega_X^n \le   C\exp (C  {\overline E}_t),
  \end{equation}
  where $A_s: =  \frac{c_t^n}{V_t} \int_{\Omega_s} (-\varphi_t - s) e^{nF_t} \omega_{X}^n$ is the energy of $(\varphi_t + s)_-$.
  \end{lemma}
\noindent{\em Proof.}
We choose a sequence of smooth positive functions $\tau_k: {\mathbb R}\to{{\mathbb R}}_+$ such that 
\begin{equation}\label{eqn:tau k}\tau_k(x) = x + \frac 1 k,\quad \mbox{ when }x\ge 0,
\end{equation}
 and $$ \tau_k(x) = \frac 1 {2k},\quad  \mbox { when } x\le -\frac {1}{k},$$ and $\tau_k(x)$ lies between $1/2k$ and $1/k$ for $x\in [-1/k, 0]$. Clearly $\tau_k$ converge pointwise to $\tau_\infty(x) = x \cdot \chi_{{\mathbb R}_+}(x)$ as $k\to\infty$, where $\chi_{{\mathbb R}_+}$ denotes the characteristic function of ${\mathbb R}_+$.
\iffalse
Consider a sequence of smooth decreasing and positive functions $f_k$ which converge pointwise to $\chi_{\Omega_s}$ and we may choose $f_k$ such that $f_k \equiv 1$ on $\Omega_s$.\fi 

We solve an auxiliary complex Monge-Amp\`ere equation on $X$
   \begin{equation}\label{eqn:aux MA}
         (\omega_t + \ddbar \psi_{t, k})^n = \frac{\tau_k(-\varphi_t - s)}{ A_{s,k}} f(\lambda[h_{\varphi_t}])^n \omega_X^n = \frac{\tau_k(-         \varphi_t - s)}{ A_{s,k}} c_t^n e^{nF_t} \omega_X^n,  
    \end{equation}
 with $\sup\psi_{t, k} = 0$
where $A_{s,k}: = \frac{c_t^n }{V_t}\int_X \tau_k(-\varphi_t - s) e^{nF_t} \omega_{X}^n$ is chosen so that the integrals of both sides of \eqref{eqn:aux MA} are equal. Note that \eqref{eqn:aux MA} admits a unique smooth solution by Yau's theorem \cite{Y}. We also observe that  as $k\to\infty$
\begin{equation}\label{eqn:limit}
A_{s,k} \to A_s =  \frac{c_t^n}{V_t} \int_{\Omega_s} (-\varphi_t - s) e^{nF_t} \omega_{X}^n, 
\end{equation}
which follows from Lebesgue's dominated convergence theorem. The limit  $A_s$ satisfies $A_s\le {\overline E}_t$, the assumed upper bound of $E_t(\varphi_t)$.

\medskip

Denote $\Phi$ to be the smooth function 
\bea
\Phi := -\varepsilon ( -\psi_{t,k} + \Lambda )^{\frac{n}{n+1}}  - (\varphi_t + s)  \eea
where
    \begin{equation}\label{eqn:choice}
0< \varepsilon: = (\frac{n+1}{n^2})^{\frac{n}{n+1}} A_{s,k}^{\frac{1}{n+1}}\gamma^{{-\frac{1}{n+1}}}, \quad 0<\Lambda := \frac{1}{(n+1) n^{n-1}} \frac{A_{s,k}}{\gamma}
     \end{equation} 
where $\gamma>0$ is the constant in the structure condition \eqref{eqn:strf} of $f$.

Since $X$ is compact without boundary, the maximum of $\Phi$ must be attained at some point, say, $x_0\in X$. If $x_0 \in X\backslash \Omega_{s}^\circ$, then 
   $$\sup_X \Phi = \Phi(x_0) =-\varepsilon ( -\psi_{t,k}(x_0) + \Lambda )^{\frac{n}{n+1}}  - (\varphi_t(x_0) + s)< -\varphi_t(x_0) - s \le 0.$$ 
 Otherwise $x_0\in \Omega_{s}^\circ$. Then at $x_0$, $\ddbar \Phi(x_0)\le 0$, and on the right hand side of \eqref{eqn:aux MA} we have $$\tau_k(-\varphi_t - s) (x_0) = - (\varphi_t(x_0) + s) + 1/k >0$$ by the definition of $\tau_k$ in \eqref{eqn:tau k}. 
 
 \smallskip
 
 We denote by $G^{i\bar j} = \frac{\partial \log f(\lambda[h])}{\partial h_{ij}} = \frac{1}{f} \frac{\partial f(\lambda[h])}{\partial h_{ij}}$ the coefficients of the linearization of the operator $\log f(\lambda[h])$ with $h =\omega_X^{-1}\cdot \omega_{t,\varphi_t}$. By the ellipticity assumption of $f(\lambda[\cdot])$, $(G^{i\bar j})$ is positive definite. Moreover, by the structure condition \eqref{eqn:strf} on $f$, we have 
 $$\det\, G^{i\bar j} = f^{-n } \det\big(\frac{\partial f(\lambda[h])}{\partial h_{ij}} \big)\ge \frac{\gamma}{f(\lambda)^n} .$$ 
Recall that the eigenvalues of $h$ are by definition $\lambda = (\lambda_1,\ldots, \lambda_n)$. Working in a basis where $h$ is diagonal and $\o_X$ the identity,
we find using the definition of $G^{i\bar j}$ that 
$$\sum_{i,j}G^{i\bar j} (\omega_{t,\varphi_t})_{\bar j i}  = \frac{1}{f(\lambda)} \sum_j \frac{\partial f(\lambda)}{\partial \lambda_j} \lambda_j= 1$$
where we have used the assumption that $f$ is homogeneous of degree one, so that $\sum_i \lambda_i\frac{\partial f(\lambda)}{\partial \lambda_i} = f(\lambda)$. 
 At the maximum point $x_0$ of $\Phi$, we can write 
\bea
0&\ge &\nonumber  G^{i\bar j}\Phi_{\bar j i}(x_0)\\
& =&\nonumber \frac{n \varepsilon}{n+1} (-\psi_{t,k} + \Lambda)^{-\frac 1{n+1}} G^{i\bar j} (\psi_{t,k})_{\bar j i} + \frac{n\varepsilon}{(n+1)^2} (-\psi_{t,k} + \Lambda)^{-\frac{n+2}{n+1}} G^{i\bar j} (\psi_{t,k})_{\bar j} (\psi_{t,k})_i - G^{i\bar j} (\varphi_t)_{\bar ji}  \\
& \ge&\nonumber  \frac{n \varepsilon}{n+1} (-\psi_{t,k} + \Lambda)^{-\frac 1{n+1}} G^{i\bar j} (\omega_ {t,\psi_{t,k}})_{\bar j i}  - G^{i\bar j} (\omega_{t, \varphi_t})_{\bar ji} + \big( 1- \frac{\varepsilon n}{n+1} (-\psi_{t, k} + \Lambda)^{-\frac{1}{n+1}}  \big) G^{i\bar j} (\omega_t)_{\bar j i} \\
& \ge&\nonumber\frac{ \varepsilon n}{n+1} (-\psi_{t, k} + \Lambda)^{-\frac{1}{n+1}} n \Big( \det G^{i\bar j}\cdot  \det (\omega_{t,\psi_{t,k}} )_{\bar j i}  \Big)^{1/n} - 1 + \big( 1- \frac{\varepsilon n}{n+1} (-\psi_{t, k} + \Lambda)^{-\frac{1}{n+1}}  \big)  G^{i\bar j} (\omega_t)_{\bar j i}\\
& \ge&\nonumber\frac{ \varepsilon n^2}{n+1} (-\psi_{t, k} + \Lambda)^{-\frac{1}{n+1}} \gamma^{1/n} \Big( \frac{\tau_k(-\varphi_t - s)}{A_{s,k}}  \Big)^{1/n} - 1 + \big( 1- \frac{\varepsilon n}{n+1}  \Lambda^{-\frac{1}{n+1}}  \big)  G^{i\bar j} (\omega_t)_{\bar j i}\\
& \ge &\nonumber \frac{\varepsilon n^2 \gamma^{1/n}}{n+1} (-\psi_{t, k} + \Lambda)^{-\frac{1}{n+1}} \Big(  \frac{-\varphi_t - s + 1/k  }{A_{s,k}}   \Big)^{1/n} -1 % + \big( 1- \frac{\varepsilon n}{n+1}  \Lambda^{-\frac{1}{n+1}}  \big) \tr_{\omega_{\varphi_t}} \omega_t
\eea
where in the third inequality we used the arithmetic-geometric inequality and in the last one we used the choice of $\varepsilon $ and $\Lambda$ in \eqref{eqn:choice}.  Thus at $x_0$ we have  
$$-(\varphi_t + s)(x_0)< A_{s,k} \Big(\frac{n+1}{n^2 \varepsilon \gamma^{1/n}} \Big)^n (-\psi_{t, k} (x_0) + \Lambda)^{n/(n+1)} =  \varepsilon (-\psi_{t, k} (x_0) + \Lambda)^{n/(n+1)}$$
which implies that $\Phi(x_0)\le 0$. Hence we can conclude that $\sup_X \Phi\le 0$, that is, on $X$
\begin{equation}\label{eqn:middle 1}\frac{- (\varphi_t + s)}{ A_{s, k}^{1/(n+1)}  } \le  (\frac{n+1}{n^2})^{\frac{n}{n+1}} \gamma^{-\frac {1}{n+1}} \big(-\psi_{t, k} + \frac{1}{(n+1) n^{n-1}} \frac{A_{s, k}}{\gamma} \big)^{\frac n{n+1}}\le C_n  \big(-\psi_{t, k} +  \frac{A_{s,k}}{\gamma} \big)^{\frac n{n+1}},
\end{equation}
for some constant $C_n$ depending only on $n$ and $\gamma$. Taking the $\big( \frac{n+1}{n}\big)$-th power of both sides of the previous equation, multiplying it by some small $\beta_0>0$, taking the exponential on both sides and then integrating the resulting inequality over $\Omega_s$, we obtain 
  \begin{equation}\label{eqn:1.8}
 \int_{\Omega_s} \exp\Big\{ \beta_0 \big(\frac{- (\varphi_t + s)}{ A_{s, k}^{1/(n+1)}  } \big)^{\frac{n+1}{n}}   \Big\} \omega_X^n \le   \exp( C_n \beta_0  A_{s, k}) \int_{\Omega_s} \exp( - C_n \beta_0 \psi_{t, k}) \omega_X^n.
  \end{equation}
 
Recall that $\omega_t + \ddbar \psi_{t, k}>0$ and $\omega_t = \chi + t \omega_X$. We may assume $\chi\le (a_0-1)\omega_X$ for some $a_0=a_0(\chi, \omega)>1$, so $\psi_{t, k}$ is also $(a_0\omega_X)$-plurisubharmonic. Now it is a basic fact in K\"ahler geometry that, for any K\"ahler class $\hat\chi$
on $X$, there is a constant $\alpha=\alpha(X,\hat\chi)$ so that
\bea
\label{eqn:alpha}
\int_X e^{-\alpha_0 \psi}\omega_X^n \le C(\alpha_0, n,\hat\chi,\omega_X)
\eea
for any $\alpha_0<\alpha$ and any $\hat\chi$-plurisubharmonic function $\psi$ with $\sup_X \psi = 0$. The local version of this statement is in \cite{H}, and the above global version in \cite{Ti}. We apply this statement with $\hat\chi=
a_0\o_X$, and fix $\alpha_0$ with $0<\alpha_0<\alpha(X,\hat\chi)$. 
 Then we choose $\beta_0 = \beta_0(n,\omega_X,\chi, \gamma)>0$ in \eqref{eqn:1.8} such that $\beta_0 C_n = \alpha_0$, and from \eqref{eqn:1.8} we can then deduce that
  \begin{equation}\label{eqn:1.9}
 \int_{\Omega_s} \exp\Big\{{  \beta_0 \big(\frac{- (\varphi_t + s)}{ A_{s, k}^{1/(n+1)}  } \big)^{\frac{n+1}{n}}   }\Big\} \omega_X^n \le   Ce^{ C A_{s, k}},
  \end{equation} 
  for some constant $C = C(n,\omega_X, \chi, \gamma)>0$.
Letting $k\to \infty$ in \eqref{eqn:1.9} we obtain from \eqref{eqn:limit} 
  \begin{equation}\label{eqn:1.10}
 \int_{\Omega_s} \exp \Big\{{  \beta_0 \big(\frac{- (\varphi_t + s)}{ A_{s}^{1/(n+1)}  } \big)^{\frac{n+1}{n}}   } \Big\} \omega_X^n \le   C e^{ C  A_{s}}\le C e^{ C{\overline E}_t}, 
  \end{equation}
   for some constant $C = C(n,\omega_X, \chi, \gamma)>0$.
The proof of Lemma \ref{lemma trudinger} is complete.
  
  \bigskip
  
We come now to the proof of Theorem \ref{thm:main} proper.
Fix $p>n$, and define $\eta: {\mathbb R}_+\to {\mathbb R}_+$ by $\eta(x) = (\log(1+x))^p$. Note that $\eta$ is a strictly increasing function with $\eta(0) = 0$, and let $\eta^{-1}$ be its inverse function. If we let
\bea
v: = \frac{\beta_0}{2}\big( \frac{-\varphi_t - s}{A_s^{1/(n+1)}} \big)^{(n+1)/n}
\eea
then we have for any $z\in\Omega_s$, by
the generalized Young's inequality with respect to $\eta$, 
   \bea
  v(z)^p e^{nF_t(z)} & \le &\nonumber \int_0^{\exp({nF_t(z)})} \eta(x) dx + \int_0^{v(z)^p} \eta^{-1}(y) dy\\
  &\le &\nonumber \exp({nF_t(z)})  (\log(1+  \exp({nF_t(z)})))^p + \int_0^{\exp({v(z) - 1})} x \eta'(x) dx\\
  &\le&\nonumber \exp({nF_t(z)}) (1+ n |F_t(z)|)^p + v(z)^p \exp({v(z)})\\
    &\le&\nonumber \exp({nF_t(z)}) (1+ n |F_t(z)|)^p + C(p) \exp({2 v(z)})
    \eea
We integrate both sides in the inequality above over $z\in \Omega_s$, and get by Lemma \ref{lemma trudinger} that %\footnote{In the inequality in fact we can get  $  \int_{\Omega_s} v^p e^{nF_t} \omega_X^n \le  2 \| e^{nF_t} \|_{L^1(\log L)^p} + \| e^{nF_t}\|_{L^{1/k}}^{1/k} + C e^{C {\overline E}_t} $ for any $k\ge 1$, by the simple calculus inequality that $e^{nF_t}\le e^{nF_t} |F_t|^p + C e^{nF_t /k}$, and this is related to the normalization of $F_t$. Recall we normalize $F_t$ such that $\int_X e^{nF_t} \omega_X^n = V$.}
   \bea
  \int_{\Omega_s} v(z)^p e^{nF_t (z)} \omega_X^n & \le&\nonumber \int_{\Omega_s} e^{nF_t } (1+ n |F_t(z)|)^p \omega_X^n + \int_{\Omega_s} e^{2v(z)} \omega_X^n\\
  & \le&\nonumber \| e^{nF_t} \|_{L^1(\log L)^p} + C + C e^{C {\overline E}_t},
  \eea
  where the constant $C>0$ depends only on $n,\omega_X, \chi, \gamma, p$. In view of the definition of $v$, this implies  
   \begin{equation}\label{eqn:final 1} \int_{\Omega_s} (-\varphi_t - s)^{\frac{(n+1)p}{n}} e^{nF_t (z) } \omega_X^n \le 2^p\beta_0^{-p} A_s^{\frac{p}{n}} \big(\| e^{nF_t} \|_{L^1(\log L)^p} + C + C e^{C  {\overline E}_t}\big).
   \end{equation}
From the definition of $A_s$ in \eqref{eqn:limit}, it follows from H\"older inequality that
\bea
     A_s  & = &\nonumber\frac{c_t^n}{V_t} \int_{\Omega_s} (-\varphi_t - s) e^{nF_t} \omega_X^n \\
     & \le&\nonumber  \Big( \frac{c_t^n}{V_t} \int_{\Omega_s} (-\varphi_t - s) ^{\frac{(n+1)p}{n}} e^{nF_t}\omega^n_X\Big)^{\frac{n}{(n+1)p}} \cdot \Big(\frac{c_t^n}{V_t}  \int_{\Omega_s} e^{nF_t} \omega_X^n \Big)^{1/q}\\
     &\le&\nonumber  A_s^{\frac {1}{n+1}}\Big( \frac{c_t^n}{V_t} 2^p\beta_0^{-p} \big(\| e^{nF_t} \|_{L^1(\log L)^p} + C + C e^{C {\overline E}_t}\big)\Big)^{\frac{n}{(n+1)p}} \cdot \Big(\frac{c_t^n}{V_t}  \int_{\Omega_s} e^{nF_t} \omega_X^n \Big)^{1/q}
\eea
where $q>1$ satisfies  $\frac{n}{p(n+1)} + \frac{1}{q} = 1$, i.e. $q = \frac{p(n+1)}{p(n+1) - n}$. The inequality above yields 
\begin{equation}\label{eqn:1.11}
   A_s\le \Big( \frac{c_t^n}{V_t} 2^p\beta_0^{-p} \big(\| e^{nF_t} \|_{L^1(\log L)^p} + C + C e^{C{\overline E}_t}\big)\Big)^{1/p} \cdot \Big(\frac{c_t^n}{V_t}  \int_{\Omega_s} e^{nF_t} \omega_X^n \Big)^{\frac{1+n}{qn}}.
\end{equation}
Observe that the exponent of the integral on the right hand of \eqref{eqn:1.11} satisfies
$$\frac{1+n}{qn} = \frac{pn + p - n}{pn} = 1+ \delta_0>1,$$
for $\delta_0: = \frac{p-n}{pn}>0$. For notation convenience, set 
\begin{equation}\label{eqn:b0}B_0 :=  \Big( \frac{c_t^n}{V_t} 2^p\beta_0^{-p} \big(\| e^{nF_t} \|_{L^1(\log L)^p} + C + C e^{C  {\overline E}_t}\big)\Big)^{1/p}.
\end{equation}
From \eqref{eqn:1.11} we then get 
\begin{equation}\label{eqn:1.14}
     A_s \le B_0 \Big(\frac{c_t^n}{V_t}  \int_{\Omega_s} e^{nF_t} \omega_X^n \Big)^{1+\delta_0}.
   \end{equation}
For any $r\in [0,1]$, we note that $-\varphi_t - s \ge r$ on $\Omega_{s+r} = \{\varphi_t \le -s -r\}$. Thus 
\begin{equation}\label{eqn:1.15}
    A_s = \frac{c_t^n}{V_t} \int_{\Omega_s} (-\varphi_t - s) e^{nF_t} \omega_X^n \ge r \cdot \frac{c_t^n}{V_t} \int_{\Omega_{s+r}} e^{nF_t} \omega_X^n.
\end{equation}
If we define $\phi:{\mathbb R}_+ \to {\mathbb R}_+$ by $$\phi(s) := \frac{c_t^n}{V_t} \int_{\Omega_s} e^{nF_t} \omega_X^n$$ then \eqref{eqn:1.14} and \eqref{eqn:1.15} imply that 
\begin{equation}\label{eqn:final 10}
r \phi(s+r) \le B_0 \phi(s)^{1+\delta_0},\quad \forall r\in [0,1] \mbox{ and } s\ge 0.
\end{equation}
$\phi$ is clearly nonincreasing and continuous, so the lemma below applies to $\phi$. It is a classic lemma due to De Giorgi, which was also used in \cite{K, EGZ}. We include a sketch of the proof for the readers' convenience, and to exhibit the dependence of $\| \varphi_t\|_{L^\infty}$ on the given data.

\begin{lemma}\label{lemma EGZ}
Let $\phi: {\mathbb R}_+\to {\mathbb R}_+$ be a decreasing right-continuous function with $\lim_{s\to\infty}\phi(s) = 0$. Assume that $r \phi(s+r)\le B_0 \phi(s)^{1+\delta_0}$ for 
some constant $B_0>0$ and all $s> 0$ and $r\in [0,1]$. Then there exists some $S_\infty = S_\infty(\delta_0, B_0, \phi)>0$ such that $\phi(s) = 0$ for all $s\ge S_\infty$.
\end{lemma}
\noindent{\em Proof.} Fix an $s_0>0$ such that $\phi(s_0)^{\delta_0} < \frac{1}{2B_0}$. This $s_0$ exists since $\phi(s)\to 0$ as $s\to\infty$. Define an increasing sequence $(s_j)$ of positive real numbers inductively by
$$s_{j+1}: = \sup\{ s> s_j |~ \phi(s) > \frac 1 2 \phi(s_j)  \}.$$
If at some stage $\phi(s_j) = 0$, we stop there. By the right-continuity of $\phi$, it follows that $\phi(s_{j+1})\le \frac{\phi(s_j)}{2}$ and $s_{j+1}\le 1 + s_j$ since $\phi(1+s_j)\le \frac{1}{2}\phi(s_j)$. It follows from the assumptions on $\phi$ that
$s_{j+1} - s_j\le 2^{-j \delta_0}$, which implies that the sequence $(s_j)$ converges to 
$$S_\infty = s_0 + \sum_{j\ge 0} (s_{j+1} - s_j)\le s_0 + \frac{1}{1-2^{-\delta_0}}.$$ It is then clear that  $\phi(s) = 0$ for any $s> S_\infty$. The lemma is proved.

\medskip

We return now to the proof of Theorem \ref{thm:main}. By Chebyshev's inequality, we have
$$\phi(s) = \frac{c_t^n}{V_t} \int_{\Omega_s} e^{nF_t} \omega_X^n \le \frac{1}{s} \frac{c_t^n}{V_t} \int_{\Omega_s} (-\varphi_t )e^{nF_t} \omega_X^n \le \frac{{\overline E}_t}{s} \to 0 \mbox{ as }s\to \infty. $$ Thus we may choose $s_0 = (2B_0)^{1/\delta_0} {\overline E}_t$ in the proof of Lemma \ref{lemma EGZ}. By \eqref{eqn:final 10} and Lemma \ref{lemma EGZ}, we deduce that 
$$\Omega_{S_\infty} = \{\varphi_t \le -S_\infty\} = \emptyset,$$
so hence
\begin{equation}\label{eqn:final}
\inf_X\varphi_t \ge -S_\infty = - (2B_0)^{1/\delta_0} {\overline E}_t - \frac{1}{1-2^{-\delta_0}},
\end{equation}
where $B_0$ is the constant in \eqref{eqn:b0} and $\delta_0: = \frac{p-n}{pn}>0$ depends only on $n$ and $p$. The proof of Theorem \ref{thm:main} is complete.

%%%%%%%%%%%%%% complete the proof of Theorem 2.
%%%%%%%%%%%%%%%%
%%%% proof of theorem 1.

\section{Proof of Theorem \ref{thm:main1}}
\setcounter{equation}{0}
\newcommand{\ent}{{\mathrm{Ent}}}
Let $\varphi$ solve the equation \eqref{eqn:main1}, where $F$ is a given smooth function and $f(\lambda[\cdot])$ is the nonlinear operator as introduced in Section 1. The setting of Theorem \ref{thm:main1} with a fixed background $\o_X$ can be viewed as a special case of the setting of Theorem \ref{thm:main} with $\chi=0$, and $t$ taken to be $1$, $\o_1=\o_X$, and in the notations (\ref{eqn:main}) and (\ref{eqn:main1}) for the two settings,
\bea
e^F=c_1e^{F_1},
\quad
V=V_1.
\eea
We observe that, in view of the normalization ${1\over V}\int_Xe^{nF_1}\o_X^n=1$ for $F_1$, 
\bea
c_1^n={1\over V}\int_Xe^{nF}\o_X^n
\leq \Ent_p(F)+e^n
\eea
for any $p\geq 1$. Thus, applying Theorem \ref{thm:main} and assuming that $\Ent_p$ is bounded, we find that $L^\infty$ bounds for $\varphi$ would follow if we can 
control the energy $E=E_{t=1}$. However, an easy application of H\"older's inequality gives
\begin{equation}
\label{eqn:energy bound}
E(\varphi)=\frac 1 V \int_X (-\varphi) e^{nF} \omega_X^n \le \Big(\frac 1 V \int_X (-\varphi)^{\frac{n}{n-1}} e^{nF} \omega_X^n \Big)^{(n-1)/n}\le C
\end{equation}
so it suffices to control the right hand side. This is done in the following theorem, part (c), which completes the proof of Theorem \ref{thm:main1}:

\begin{theorem}
\label{lemma entropy}\label{thm:energy}
If $f(\lambda[\cdot])$ satisfies the structure condition \eqref{eqn:strf}, then the following holds:

\medskip

{\rm(a)} Assume that $p\in (0,n)$. Then there exist constants $c_p$, $C_p>0$ depending only on  $\o_X$, $n$, $p$, $\gamma$ and the generalized entropy $\ent_p(F)$ such that
\bea \label{eqn:main a}
\int_X {\rm exp}\big\{c_p(-\varphi)^{n\over n-p}\big\}\o_X^n
\leq C_p.
\eea

{\rm (b)} Assume that $p = n$. Then for any $N>0$, there exists constants $c_{N}>0$, $C_{N}>0$ depending on $n,\o_X$, $N$, $\gamma$, and the generalized entropy $\ent_n(F)$ so that
\bea\label{eqn:main b}
\int_X{\rm exp}\big\{c_{N}(-\varphi)^N\big\}\o_X^n\leq C_{N}.
\eea

{\rm (c)} We have the energy estimate:
\bea\label{thm c}
\int_X(-\varphi)^N e^{nF}\o_X^n \leq C
\eea
for $N = {n\over n-p}$ if $p\in [1,n)$, and for 
any $N>0$ if $p= n$, where the constant $C$ on the right hand side of \eqref{thm c} depends on $n, \omega_X, \gamma, N$ and the entropy $\ent_p(F)$.
\end{theorem}

\bigskip

We observe that, in the special case of the Monge-Amp\`ere equation and when $p=1$, these estimates have been established in \cite{DGL}, using pluripotential theory.  Further if $p>n$, then Theorem \ref{thm:main1} implies the solutions are $L^\infty$ bounded. Theorem  \ref{thm:energy} gives a more complete integral estimate of such solutions for the full range of $p\in (0,\infty)$, using pure PDE methods.

For the proof of Theorem \ref{thm:main1}, we only need part (c) of Theorem \ref{lemma entropy}, but we give the proofs of the other parts as well, as they are of independent interest. As mentioned in the Introduction, for this we need the ABP method developed in \cite{B11} and \cite{CC} in order to get better integral bounds for $\varphi$ as stated in \eqref{eqn:main a} and \eqref{eqn:main b}. In particular, we have fix a background metric $\omega_X$, as this method is not effective for handling degenerating families.

\bigskip

%\noindent{\em Proof of Lemma \ref{lemma entropy}.}
\newcommand{\innpro}[1]{\langle #1 \rangle}
%%%%%%%%%%%%%%%%%%%%%%%%%%
%%%%%%%%%%%%%%%%%%%%%%%%%%

\medskip

%%%%%%%%%%%%%%%%%%%%%%%%%%%
%%%%%%%%%%%%%%%%%%%%%%%%%%%%

Suppose $p\in (0, n]$ and write $$\Psi_p: = \frac 1 V \int_X (F^2+1)^{p/2} e^{nF}\omega_X^n \sim\frac 1 V \int_X |F|^p e^{nF} \omega_X^n = \ent_p(F). $$

\smallskip

\noindent We  can solve the complex Monge-Amp\`ere equation
\begin{equation}\label{eqn:MA 2}(\omega_X +\ddbar \psi)^n = \frac{(F^2 +1)^{p/2}}{\Psi_p} e^{nF}\omega_X^n,\quad \sup_X\psi = 0.\end{equation}

\begin{lemma}\label{lemma 2} Let $\varphi, \psi$ be the solutions to the equations \eqref{eqn:main1}, \eqref{eqn:MA 2}, respectively. There exist constants $\Lambda, \lambda, C>0$ depending on $n, \omega_X, \gamma, p, $ and $\Psi_p$  such that
$$\sup_X ( - (-\psi +\Lambda)^\beta - \lambda \varphi ) \le C$$
where  $\beta = \frac{n-p}{n}\in (0,1)$ if $p<n$, and if $p=n$, we choose an arbitrary $\beta = N^{-1} \in (0,1)$ and the constants $\lambda, \Lambda, C$ depend additionally on $N$.
\end{lemma}

 In the proof below, for a smooth function $u$ on $X$, we denote \begin{equation}\label{eqn:G} \Box u := G^{i\bar j} u_{\bar j i}, \quad |\nabla u|_G^2 := G^{i\bar j} \nabla_{\bar j} u\nabla_i u\end{equation} and $$\tr_G \alpha: = G^{i\bar j} \alpha_{\bar j i},\quad \mbox{ for a smooth $(1,1)$-form $\alpha = \alpha_{\bar j i}\sqrt{-1} dz^i\wedge d\bar z^j$}, $$ 
where as in Lemma \ref{lemma trudinger}, $G^{i\bar j} = \frac{\partial \log f(\lambda[h])}{\partial h_{ij}}$ is the coefficient matrix of the linearized operator of $\log f(\lambda[\cdot])$ at  $h = \omega_X^{-1}\cdot \omega_\varphi$, and by the structure condition \eqref{eqn:strf} on $f$, we have $\det (G^{i\bar j}) \ge \gamma f^{-n}$ and $G^{i\bar j}$ is positive definite.

\bigskip

\noindent{\em Proof of Lemma \ref{lemma 2}.} We choose constants as follows:
\begin{equation} \label{eqn:constants chosen}
\Lambda  = \Big( \frac{4^n}{n^n \gamma } \frac{2^p}{\alpha_0^p} \Psi_p  \Big)^{1/n(1-\beta)}   ,\quad  \lambda = 4\beta \Lambda^{-(1-\beta)}
\end{equation}
where as usual $\alpha_0=\alpha_0(X,\omega_X)$ is a fixed constant smaller than the $\alpha$-invariant of $(X,\omega_X)$. 
 We denote 
\begin{equation}\label{eqn:beta}\rho := -(-\psi + \Lambda)^\beta -\lambda \varphi .\end{equation} %and $H = e^\rho$.

For notation convenience, we set $\phi_\delta(t) = t + \sqrt{t^2 + \delta}>0$, which is a smoothing of $2 \max(t, 0)$ and converges to it as $\delta\to 0$. Here we will first fix a $\delta>0$ small and later on $\delta$ will be sent to zero. All constants appearing in the proof are independent of $\delta$, unless stated otherwise. From now on we will consider $\phi_\delta(\rho)$ which is monotone decreasing and converges to $2 \rho_+$ as $\delta\to 0$. We will omit the $\delta$ in $\phi_\delta(\rho)$ and simply write $\phi(\rho)$.

\medskip

\noindent We define a {\em smooth} function 
$$H = \phi(\rho)^b,$$
where  $b = 1 + \frac 1{4n}>1$ is constant. Since $X$ is compact, $H$ must achieve its maximum at some point in $X$, say, $x_0 $, and we denote $\sup_X H =: M>0$ (if $M=0$ there is nothing to prove). Let $r = \min\{1, r(X,\omega_X)\}$ where $r(X,\omega_X)>0$ is the injectivity radius of $(X,\omega_X)$ viewed as a compact Riemannian manifold. So we can identify the geodesic ball $B_r(x_0)$ as an open smooth domain in ${\mathbb R}^{2n}$ with Euclidean diameter bounded by $3r$, say. Let $\theta \in (0,1)$ be a small constant defined by
\begin{equation}\label{eqn:theta}
\theta : = \min\{ \frac{r^2\beta\Lambda^{-(1-\beta)}}{100 M^{1/b}}, \frac{r^2}{100n}\}<\frac 1{10}.
\end{equation} 
As in \cite{CC}, we choose an auxiliary smooth function $\eta$ defined on $B_r(x_0)$ so that $\eta \equiv 1$ on $B_{r/2}(x_0)$ and $ \eta \equiv 1 -\theta$ on $X\backslash B_{3r/4}(x_0)$, and $\eta$ lies between $1$ and $1-\theta$ in the annulus $B_{3r/4}(x_0)\backslash B_{r/2}(x_0)$. Moreover $\eta$ can be chosen to  satisfy 
\bea \label{eqn:eta 1}|\nabla \eta|_g^2 \le \frac{10 \theta^2}{r^2},\quad |\nabla^2 \eta|_g \le \frac{10\theta}{r^2},\eea
where we identify $\omega_X$ with its associated Riemannian metric $g$.

\medskip

\newcommand{\fpp}{\varepsilon \psi - \lambda \varphi}

We can now calculate,
\begin{equation}\label{eqn:1}
 \Box( H \eta ) = \eta  \Box H  + H  \Box \eta + 2 Re\big( G^{i\bar j} \nabla_{\bar j} H \nabla_i \eta   \big).\end{equation}
We observe that the middle term in \eqref{eqn:1} satisfies % (below for a function $u$, $|\nabla u|_G^2 := G^{i\bar j} \nabla_{\bar j} u\nabla_i u$ and $\tr_G \alpha: = G^{i\bar j} \alpha_{\bar j i}$ for a smooth $(1,1)$-form $\alpha$)
$$H  \Box \eta = H \tr_G \ddbar \eta\ge - H \frac{10\theta}{r^2} \tr_G \omega_X.$$
The last term in \eqref{eqn:1} satisfies
\bea
2 Re\big( G^{i\bar j} \nabla_{\bar j} H \nabla_i \eta   \big) & = & \nonumber 2b \phi(\rho)^{b-1} Re \big( G^{i\bar j}\nabla_{\bar j} \phi(\rho) \nabla_i \eta   \big)\\
& \ge & \nonumber -\frac{b (b-1)}{2} \phi(\rho)^{b-2}| \nabla \phi(\rho)  |^2_G  - \frac{2b}{b-1}\phi(\rho)^b  |\nabla \eta|^2 _G\\
& \ge& \nonumber -\frac{b (b-1)}{2} \phi(\rho)^{b-2}| \nabla \phi(\rho)  |^2_G - \frac{2b}{b-1} \phi(\rho)^b  \frac{10\theta^2}{r^2} \tr_G \omega_X
\eea 
where in the first inequality we applied the Cauchy-Schwarz inequality.
The first term in \eqref{eqn:1} is
\bea\label{eqn:middle 12 new}
\eta \Box H &=&b\eta\phi(\rho)^{b-1}  \Box \phi(\rho) + b (b-1) \eta \phi(\rho)^{b-2} | \nabla \phi(\rho)|_G^2\\
& = & 
b\eta\phi(\rho)^{b-1}  \phi'(\rho) \Box \rho + b \eta \phi(\rho)^{b-1} \phi''(\rho) |\nabla \rho |^2_G    + b (b-1) \eta \phi(\rho)^{b-2} | \nabla \phi(\rho)|_G^2. 
\nonumber
\eea
We note that the middle term in \eqref{eqn:middle 12 new} is nonnegative due to the fact that $$\phi''(t) = \frac{1}{\sqrt{t^2 + \delta}} - \frac{t^2}{(\sqrt{t^2 + \delta})^3} = \frac{\delta}{(t^2 + \delta)^{3/2}}>0. $$
To deal with the first term in \eqref{eqn:middle 12 new} we note by the homogeneity of degree one assumption on $f$ that
$$\Box \varphi = \tr_G \ddbar \varphi = \tr_G \omega_\varphi - \tr_G \omega_X = 1 - \tr_G \omega_X.$$
Then we  
calculate
%
%\smallskip
%
\bea \label{eqn:mid 10}
\Box \rho & =\nonumber & \Box( - (-\psi + \Lambda)^\beta -\lambda \varphi)\\
 & =\nonumber &  \beta (-\psi + \Lambda)^{\beta -1} \Box \psi +  \beta (1-\beta) (-\psi + \Lambda)^{\alpha -2} |\nabla \psi|_G^2 - \lambda \Box\varphi \\
&\ge & \nonumber  \beta ( -\psi + \Lambda  )^{\beta - 1} \tr_G \omega_\psi -   \beta ( -\psi + \Lambda  )^{\beta - 1} \tr_G \omega_X -\lambda + \lambda \tr_G \omega_X\\
&\ge &\nonumber   n \gamma^{1/n}  \beta ( -\psi + \Lambda  )^{\beta - 1} \big( \frac{(F^2+1)^{p/2}}{\Psi_p}   \big)^{1/n} +(\lambda -   \beta ( -\psi + \Lambda  )^{\beta - 1}) \tr_G \omega_X -\lambda\\
&\ge &   n \gamma^{1/n}  \beta ( -\psi + \Lambda  )^{\beta - 1} \big( \frac{|F|^p}{\Psi_p}   \big)^{1/n} +(\lambda -   \beta \Lambda  ^{\beta - 1}) \tr_G \omega_X -\lambda ,
\eea
%%%%%%
where in the second inequality we used the arithmetic-geometric inequality and  the equations for $\varphi$ and $\psi$. Plugging these inequalities into \eqref{eqn:1}, we obtain 
%%%%%%%%%%
\bea\label{eqn:text 1}
\Box( H \eta ) & \ge & \nonumber  - \frac {10 \theta}{r^2} \phi(\rho)^b  \tr_G \omega_X - \frac{2 b}{b-1} \phi(\rho)^b \frac{10\theta^2}{r^2} \tr_G\omega_X \\
&\qquad & \nonumber + b \eta\phi(\rho)^{b-1} \phi'(\rho) \Big( n \gamma^{1/n}  \beta ( -\psi + \Lambda  )^{\beta - 1} \big( \frac{|F|^p}{\Psi_p}   \big)^{1/n} +(\lambda -   \beta \Lambda  ^{\beta - 1}) \tr_G \omega_X -\lambda   \Big) \\
& \ge&  b \phi(\rho)^{b-1} \Big\{\frac 9{10} (\lambda -  \beta \Lambda^{\beta - 1}) \phi'(\rho) \tr_G \omega_X - \frac{20\theta}{r^2 b}\phi(\rho) \tr_G \omega_X \\
&\qquad& \nonumber  - \frac{20\theta^2}{(b-1) r^2} \phi(\rho) \tr_G \omega_X + n \gamma^{1/n}  \beta ( -\psi + \Lambda  )^{\beta - 1} \big( \frac{|F|^p}{\Psi_p}   \big)^{1/n} - \lambda \phi'(\rho)   \Big\},%\\
%& \ge  q (\fpp)^{q-1}_+ \Big(  c(n,k)  \phi(\rho)^{1/n} A_F^{-1/n} - C   \Big),
\eea
%%%%%%%%%%%%%%%%
where we ignored $\eta$ in the last inequality since $\eta>9/10$.
To deal with the right hand side in the equation \eqref{eqn:text 1}, we note that on the set $\{\rho\le 0\}$ $$0\le \phi(\rho) = \rho + \sqrt{\rho^2 + \delta} = \frac{\delta}{\sqrt{\rho^2 + \delta} - \rho}\le \sqrt{\delta},$$
and on this same set the function $\phi'(\cdot)$ satisfies $$1\ge \phi'(\rho) = 1 + \frac{\rho}{\sqrt{\rho^2 + \delta}} = \frac{\phi(\rho)}{\sqrt{\rho^2 + \delta}}\ge 0.$$
So on the set $\{\rho\le 0\}$ the right hand side of \eqref{eqn:text 1} is greater or equal to
\bea \nonumber
b \phi(\rho)^{b-1} \Big( - \frac{20\theta}{r^2 b}\sqrt{\delta} \tr_G \omega _X - \frac{20\theta^2}{(b-1) r^2}\sqrt{\delta} \tr_G \omega_X - \lambda   \Big)
\eea
On the other hand, on the set $\{\rho>0\}$, we know $2\ge \phi'(\rho)>1$, so the first three terms on the right hand side of \eqref{eqn:text 1} are positive due to the choice of $\theta$ in \eqref{eqn:theta}, and $ \Lambda, \lambda$ in (\ref{eqn:constants chosen}) and the fact that $\phi(\rho)\le M^{1/b}$. So  on the set $\{\rho>0\}$ the right hand side of \eqref{eqn:text 1} is greater or equal to
$$b \phi(\rho)^{b-1} \Big( n \gamma^{1/n}  \beta ( -\psi + \Lambda  )^{\beta - 1} \big( \frac{|F|^p}{\Psi_p}   \big)^{1/n} - \lambda    \Big).$$
%where in the last inequality we use the choice of $\theta$ in \eqref{eqn:theta} and the fact that $\phi(\rho)\le M^{1/q}$.

Combining the above two cases, we obtain 
\bea
\Box(H\eta) & \ge & \nonumber  b \phi(\rho)^{b-1} \Big( n \gamma^{1/n}  \beta ( -\psi + \Lambda  )^{\beta - 1} \big( \frac{|F|^p}{\Psi_p}   \big)^{1/n} - \lambda    \Big)\cdot \chi_{\{\rho>0\}} \\
&\qquad & \nonumber -  b \phi(\rho)^{b-1} \Big(  \frac{20\theta}{r^2 b}\sqrt{\delta} \tr_G \omega _X  + \frac{20\theta^2}{(b-1) r^2}\sqrt{\delta} \tr_G \omega_X + \lambda   \Big)\cdot\chi_{\{\rho\le 0\}},
\eea
where $\chi_E$ denotes the characteristic function of a set $E$.

\medskip

We now apply the ABP maximum principle to the function $H\eta$ on the domain $B_r(x_0)=:B_0$ in ${\bf R}^{2n}$. It follows that
{\small
\bea \label{eqn:ABP new}
\sup_{B_0} (H\eta) & \le & \nonumber \sup_{\partial B_0} (H \eta) + C(n) r\Big\{ \int_{B_0\cap \{\rho>0\}} \frac{ \phi(\rho)^{2n(b-1)}   ( n \gamma^{1/n}  \beta ( -\psi + \Lambda  )^{\beta - 1} \big( \frac{|F|^p}{\Psi_p}   \big)^{1/n}  - \lambda    \big)_-^{2n} }{ (\det G^{i\bar j})^2  }   \omega_X^n \\
&\qquad& \nonumber + \int_{ B_0\cap \{\rho\le 0\}   }  \frac{ \phi(\rho)^{2n(b-1)}  ( \frac{20\theta}{r^2 b}\sqrt{\delta} \tr_G \omega _X  + \frac{20\theta^2}{(b-1) r^2}\sqrt{\delta} \tr_G \omega_X + \lambda   )^{2n}      }{(\det G^{i\bar j})^2}    \omega_X^n     \Big\}^{1/2n}\\
& \le& \nonumber \sup_{\partial B_0} (H \eta) + C(n) \Big\{ \int_{B_0\cap \{\rho>0\}} \frac{\sqrt{ \phi(\rho)}   \big( n \gamma^{1/n} \beta ( -\psi + \Lambda  )^{\beta - 1} \big( \frac{|F|^p}{\Psi_p}   \big)^{1/n}  - \lambda    \big)_-^{2n} }{ e^{-2nF}  }   \omega_X^n    \\
&\qquad& +C' \delta^{n(b-1)}        \Big\}^{1/2n}
\eea}%
where the constant $C'=C'(n, F, G, \omega_X)$ in the last term may not be uniformly bounded, but this is not a concern, since later on we will let $\delta\to 0$.
We observe that  the integral involved in the last inequality is in fact integrated over the set where $ n \gamma^{1/n} \beta ( -\psi + \Lambda  )^{\beta - 1} \big( \frac{|F|^p}{\Psi_p}   \big)^{1/n}  - \lambda <0$ and $\rho >0$, and over this set we have from the choice of constants in (\ref{eqn:constants chosen})
%%%
 $$|F|\le (\Psi_p)^{1/p} \Big(\frac{\lambda}{n\gamma^{1/n} \beta} \Big)^{n/p} (-\psi + \Lambda)^{(1-\beta)n/p} = \frac{\alpha_0}{2}  (-\psi + \Lambda)^{(1-\beta)n/p}.$$
 %%%%
% by the definition of $\Phi(F) = \sqrt{F^2 + 1}$. 
 At the same time, on the same set, we have $0<\rho \le   - \lambda \varphi$ and $\phi(\rho)\le 2\rho + \sqrt{\delta}$. Therefore, we obtain from (\ref{eqn:ABP new}) that
  \bea\label{eqn:C0 new}
  M &  \le & \nonumber (1-\theta) \sup_{\partial B_0} H + C\Big(\int_{B_0}  ( - \varphi + \sqrt{\delta})^{1/2} \exp\Big( \frac{\alpha_0}{2} (-\psi + \Lambda)^{(1-\beta)n/p}\Big)         \omega_X^n + C' \delta^{n(b-1)}   \Big)^{1/2n}\\ 
 &  \le & \nonumber (1-\theta) M + C\Big(\int_{B_0}  \left( - \varphi + \exp\Big( {\alpha_0}(-\psi + \Lambda)^{(1-\beta)n/p}\Big) \right)         \omega_X^n + C' \delta^{n(b-1)}   \Big)^{1/2n}\\ 
   & \le &  (1-\theta) M + C_0 +C' \delta^{(b-1)/2}
 \eea
where $C_0 = C_0(n,\omega_X,\gamma, \Psi_p,p )>0$ is independent of $\delta$ and in the last inequality we have used the following inequalities:

\medskip

(1)  $\int_X (-\varphi)\omega_X^n \le C(n,\omega_X)$ which follows from the Green's formula, $n + \Delta_{\omega_X} \varphi >0$ (since $\lambda[h_\varphi]\in \Gamma\subset \{\lambda_1+\cdots + \lambda_n>0\}$) and the normalization condition $\sup_X \varphi = 0$.

\smallskip 

(2) $\int_X \exp\Big( {\alpha_0}(-\psi + \Lambda)^{(1-\beta)n/p}\Big) \omega_X^n\le C(n, \omega_X, \Psi_p, p)$. By the choice of $\beta$, if $p<n$, $(1-\beta)n/p = 1$; and if $p=n$, then $(1-\beta)n/p = 1- N^{-1}<1$, Young's inequality gives the desired estimate.

\smallskip

Hence we conclude that with the choice of $\theta$ in \eqref{eqn:theta}
$$\min(M^{1-\frac 1 b}, M)\le C_0 + C'  \delta^{(b-1)/2}
\, \,\Rightarrow \, M \le C_0 + C' \delta^{2 b},$$ which implies that 
$$\sup_X  2 \rho_+ \le \sup_X \phi(\rho) = M \le  C_0 + C' \delta^{2 b}.$$
Finally letting $\delta\to 0$ yields the desired estimate $\sup_X \rho_+\le C_0$ for some positive constant $C_0= C_0(n,\omega_X,\gamma, \Psi_p)$ (which may be different from the $C_0$ in (\ref{eqn:C0 new})). The proof of Lemma \ref{lemma 2} is complete.

\bigskip

\noindent{\em Proof of Theorem \ref{lemma entropy}}. (a) and (b) follow easily from Lemma \ref{lemma 2} and the fact the $\omega_X$-plurisubharmonic function $\psi$ satisfies $\int_X e^{-\alpha_0\psi}\omega_X^n\le C(n,\omega_X)$.

\medskip

The inequality \eqref{thm c} in (c) is an immediate consequence of the estimates in (a) and (b), and Jensen's inequality. More precisely, we have
$$\frac  1V \int_X e^{- n F + c_0 (-\varphi)^N} e^{n F}\omega_X^n =\frac 1 V \int_X e^{c_0 (-\varphi)^N}\omega_X^n \le C(n, \omega_X, \gamma,\Psi_p, N).$$
Taking the logarithms of both sides and applying Jensen's inequality yields $$\frac{1}{V} \int_X \big( c_0 (-\varphi)^N - n F \big) e^{nF} \omega_X^n\le C(n,\omega_X,\gamma, \Psi_p, N),$$
from which \eqref{thm c} follows after noting that $\Psi_p$ is equivalent to the entropy $\ent_p$ and if $p\ge 1$
$$\frac{1}{V} \int_X F e^{nF} \omega_X^n \le \ent_p.$$
The proof of Theorem \ref{lemma entropy} is complete.

%%%%%%%%%%%%%%%%%%%%%%%%%%%%%%
%%%%%%%%%%%%%%%%%%%%%%%%%%%%%%%%

%%%%%%%%%%%%%%%%%%%%%%%%%%
%%%%%%%%%%%%%%%%%%%%%%%55555

%%%%%%%%%%%%%%%%%%%%%%%%%%%%%%%%%%%%%%%%%%%%%%%

\section{Monge-Amp\`ere equations}
\setcounter{equation}{0}
\renewcommand{\eqref}[1]{\ref{#1}}

In this and the next section, we apply Theorems \ref{thm:main1} and \ref{thm:main} to the specific cases of the Monge-Amp\`ere and Hessian equations on a compact K\"ahler manifold $(X,\o_X)$.

\smallskip
We begin by noting that the structural condition (\ref{eqn:strf}) holds for many equations and is usually easy to check:

\begin{lemma}
\label{lemma structure}
Assume that $f:{\mathbb R}^n \to {\mathbb R}_+$ is a concave and homogeneous function of degree one, which satisfies $\frac{\partial f(\lambda)}{\partial \lambda_j}>0$ for any $\lambda$ in an admissible cone $\Gamma\subset {\mathbb R}^n$. Assume that there is a $\gamma>0$ such that 
\begin{equation}\label{eqn:homo}
f(\mu)\ge n\gamma^{1/n} (\prod_j \mu_j)^{1/n}, \quad \mbox{for all }\mu\in \Gamma_n:=  \{\lambda\in {\mathbb R}^n: \lambda_1>0,\ldots, \lambda_n>0\}.
\end{equation}
Then $f$ satisfies the structural condition (\ref{eqn:strf}).
\end{lemma}
\noindent{\em Proof.} By the concavity of $f$ on $\Gamma$, for any $\lambda,\mu\in \Gamma$ we have
\begin{equation}
\label{eqn:s1 1}f(\mu) \le f(\lambda) + \sum_{j=1}^n (-\lambda_j +\mu_j) \frac{\partial f(\lambda)}{\partial \lambda_j} = \sum_{j=1}^n\mu_j \frac{\partial f(\lambda)}{\partial \lambda_j},\end{equation}
where we have used the homogeneity of degree one assumption on $f$, which implies that $\sum_j \lambda_j \frac{\partial f(\lambda)}{\partial \lambda_j} = f(\lambda)$. Taking the infimum of the right hand side of (\ref{eqn:s1 1}) over all $\mu\in \Gamma_n$ with $\prod_{j=1}^n \mu_j = 1$, by the arithmetic-geometric inequality  and the assumption (\eqref{eqn:homo}) on $f$, we get
$$\prod_{j=1}^n \frac{\partial f(\lambda)}{\partial \lambda_j}\ge n^{-n} \big\{\inf_{\mu\in \Gamma_n, \prod_j \mu_j = 1} f(\mu)\big\}^{n} \ge \gamma>0.
$$
The desired inequality (\ref{eqn:strf}) follows from this inequality upon diagonalizing the matrix $h$. The proof of Lemma \ref{lemma structure} is complete.

\medskip
It follows immediately from this lemma that the functions $f(\lambda)=(\prod_{j=1}^n\lambda_j)^{1\over n}$, $f(\lambda)=\sigma_k(\lambda)^{1\over k}$, and $f(\lambda)
=({\sigma_k(\lambda)\over\sigma_\ell(\lambda)})^{1\over k-\ell}+c\sigma_p(\lambda)^{\frac 1 p}$ for $c>0$, $n\geq k\geq\ell\geq 1$, $n\geq p\geq 1$, corresponding respectively to the Monge-Amp\`ere, the Hessian, and the quotient Hessian equations, all satisfy the structural condition (\eqref{eqn:strf}), and the constant $\gamma$ depends only on the given numbers $n, k, \ell, p, c>0$, and the admissible cone $\Gamma = \Gamma_k=\{\lambda\in {\mathbb R}^n: \sigma_1(\lambda)>0,\ldots, \sigma_k(\lambda)>0\}$. The last equation appeared recently in \cite{CSz}.

\smallskip
In this section, we focus on the Monge-Amp\`ere equation. To discuss the underlying geometry, it is convenient to rewrite it in the more usual form
\bea
\label{eqn:MA}
(\o_t+i\p\bar\p\varphi_t)^n=c_t^n e^{nF_t}\o_X^n,
\quad \lambda[h_{t,\varphi_t}]\in\Gamma,
\quad \sup_X\varphi_t=0,
\eea
where $\o_t = t\omega_X+\chi$ is the family of degenerating background metrics. To apply Theorem \ref{thm:main}, we need to control the ratio $c_t^n/V_t$ and the energy $E_t$. This 
can be readily done using the following two easy lemmas:

\begin{lemma}\label{lemma volume}
Let $V= \int_X \omega_X^n$ be the volume of $(X,\omega_X)$, then 
$$V^{-1} =  \frac{c_t^n}{V_t},\quad \forall t\in (0,1].$$
\end{lemma}
\noindent{\em Proof.} Integrating both sides of (\ref{eqn:MA}), we get
$$c_t^n V = c_t^n \int_X e^{nF_t} \omega_X^n = \int_X (\omega_t + \ddbar \varphi_t)^n = \int_X \omega_t^n = V_t.$$

\bigskip

\begin{lemma}\label{lemma 5}
There is a uniform constant $C>0$ depending only on $n, \| e^{nF_t}\|_{L^1(\log L)^1(\omega_X^n)}$, $\omega_X$, $\chi$ such that for all $t\in (0,1]$
\bea
E_t(\varphi_t)\leq C.
\eea
\end{lemma}
\noindent{\em Proof.} 
Recall that 
$$E_t(\varphi_t) = \frac{c_t^n}{V_t} \int_X (-\varphi_t) e^{nF_t}\omega_X^n = \frac{1}{V_t} \int_X (-\varphi_t) \omega_{\varphi_t}^n,
$$
so that it suffices to show that
$\frac 1 {V_t} \int_X (-\varphi_t) \omega_{\varphi_t}^n \le C$, $\forall t\in (0,1].$
We may assume without loss of generality that $\chi\le \omega_X$, so the $\omega_t$-plurisubharmonic function $\varphi_t$ is also $2\omega_X$-plurisubharmonic  and by the $\alpha$-invariant estimate, there is an $\alpha_0 = \alpha_0(X,\omega_X)>0$ such that 
$$ \frac 1{V_t}\int_X \exp\Big( -\log \frac{\omega_{\varphi_t}^n}{\omega_X^n}  -\alpha_0 \varphi_t    \Big)\omega_{\varphi_t}^n   =\frac{1}{V_t}\int_X e^{-\alpha_0 \varphi_t}\omega_X^n\le \frac{C(n,\omega_X)}{V_t}$$
By Jensen's inequality it follows that
$$ \frac{1}{V_t} \int_X\Big( -\log \frac{\omega_{\varphi_t}^n }{ \omega_X^n} - \alpha_0\varphi_t  \Big) \omega_{\varphi_t}^n \le \log  C - \log V_t   $$
which implies that 
\bea
\frac{1}{V_t} \int_X (-\alpha_0 \varphi_t)\omega_{\varphi_t}^n & \le&\nonumber \frac{1}{V_t} \int_X \log (e^{nF_t} c_t^n) \omega_{\varphi_t}^n + \log C - \log V_t \\
& =&\nonumber \int_X (nF_t) e^{nF_t} \frac{c_t^n}{V_t} \omega_X^n + \log C + \log \frac{c^n_t}{V_t}\\
&\le&\nonumber C \| e^{nF_t}\|_{L^1(\log L)^1(\omega_X)} + C,
\eea
from which the estimate follows since $c^n_t$ and $V_t$ are uniformly equivalent by Lemma \ref{lemma volume}. 

\bigskip
Theorem \ref{thm:main} together with the preceding lemmas implies at once the following basic estimates of Kolodziej \cite{K}, Eyssidieux, Guedj, and Zeriahi \cite{EGZ}, and Demailly and Pali \cite{DP}:

\begin{theorem}
\label{thm:MA}
Consider the above family (\ref{eqn:MA}) of complex Monge-Amp\`ere equations, with respect to the degenerating background metrics $\o_t$, $t\in (0,1]$. Fix any $q>n$. If $\varphi_t$ is a family of $C^2$ solution, normalized by $\sup_X\varphi_t=0$, and if $ \| e^{nF_t}\|_{L^1(\log L)^q(\omega_X^n)}$ is uniformly bounded in $t$, then $\|\varphi_t\|_{L^\infty(X)}$ is uniformly bounded in $t$ as well.
\end{theorem}

 In many applications, the $(1,1)$-form $\chi$ is chosen to be $\chi=\pi^*\o_Y$, where $\pi:X\to Y$ is a holomorphic map between K\"ahler manifolds and $\omega_Y$ is a K\"ahler metric on $Y$.

\section{Fully non-linear and Hessian equations}
\setcounter{equation}{0}

In this section, we consider applications of Theorem \ref{thm:main} to fully non-linear equations besides Monge-Amp\`ere equations. For these applications, we need to consider the relative volumes $c_t^nV_t^{-1}$ and the energies $E_t(\varphi_t)$. We begin with a simple estimate for $E_t(\varphi_t)$, which generalizes the simple considerations which applied earlier to Monge-Amp\`ere equations and is a straightforward application of the H\"older inequality,

\begin{lemma}
\label{eqn:entropy bound}
Consider the energy $E_t(\varphi_t)$ in the formalism for degenerating background metrics as in the set-up (\ref{eqn:main}). Then we have for any $p,q> 1$ and ${1\over p}+{1\over q}=1$,
\bea
E_t(\varphi_t)
\leq {c_t^n\over V_t}\,\|e^{nF_t}\|_{L^q}\,\|\varphi_t\|_{L^p}.
\eea
\end{lemma}

%This is a straightforward application of the H\"older inequality,
%\begin{equation}\label{eqn:cor need}
%E_t(\varphi_t) = \frac {c_t^n}{V_t} \int_X (-\varphi_t) e^{nF_t} \omega_X^n\le  \frac {c_t^n}{V_t}\Big( \int_X e^{pn F_t} \omega_X^n\Big)^{1/p} \Big(\int_X (-\varphi_t)^{q}\omega_X^n  \Big)^{1\over q}\end{equation}
%for some constant $C=C(n,\omega_X, \chi, q, K)$. The last inequality in (\ref{eqn:cor need}) holds because $\frac{q}{q - 1} < \frac{n}{n-k}$ by the choice of $q$.

Next we note the following uniform $L^p$ estimate for general functions $u$, whose Hessian is in a cone $\Gamma_k=\{\lambda; \sigma_\ell>0,\ 1\leq\ell\leq k\}$, which is an analogue of the $\alpha$-invariant estimate for plurisubharmonic functions. It is well-known to experts, but we supply a statement and proof without pluripotential theory, as we could not find a convenient reference:

\begin{lemma}\label{lemma 6}\label{lemma new 8}
 For any $p\in (0, \frac{n}{n-k})$, there is a uniform constant $C=C(n, p, \omega_X)>0$ such that 
$$\| u\|_{L^p(\omega_X^n)}\le C, $$ with $\sup_X u = 0$ and $\lambda[\omega_{u}]\in \Gamma_k$, $\omega_u = \omega_X+\ddbar u$.
\end{lemma}
\noindent{\em Proof.} We use an idea in \cite{DK}.  %First observe that $\lambda[h_{2 \omega_X + \ddbar \varphi_t}] \in \Gamma_k$, so we have $\varphi_t$ is $2\omega_X$ $k$-subharmonic for any $t\in (0,1]$. 
Without loss of generality we may assume ${\mathrm{Vol}}(X,\omega_X) = 1$.

\medskip

Fix an $s>0$ and a small $\epsilon> 0$. Let $K=\{u \le -s\}\subset X$ be compact sub-level set of $u$. We choose a sequence of smooth positive functions $\eta_j$ with $\int_X \eta_j \omega_X^n = 1$ which converge to $\eta_\infty:=a V_K^{2\epsilon - 1}\chi_{K} + a \cdot\chi_{X\backslash K }$ in $L^{1+\epsilon}(\omega_X^n)$ and also pointwise, where $V_K = \int_K \omega_X^n$ is the volume of the set $K$ and $a>0$ is a constant such that 
$$\int_X \eta_\infty\omega_X^n = \int_X (a V_K^{2\epsilon - 1}\chi_{K} + a \cdot\chi_{X\backslash K }) \omega_X^n = aV_K^{2\epsilon} + a {\mathrm{Vol}(X\backslash K)} = 1.$$
It is not hard to see that $1/2\le a\le \max(2,2^{2\epsilon}) = 2$. Hence 
$$\int_X \eta_\infty^{1+\epsilon} \omega_X^n = a^{1+\epsilon} V_K^{\epsilon + 2\epsilon^2} + a^{1+\epsilon} {\mathrm{Vol}(X\backslash K)} \le 4. $$ Thus we may assume $\| \eta_j \|_{L^{1+\epsilon}(\omega_X)} \le 5$ for $j$ large enough. We solve the complex Monge-Amp\`ere equations
$$(\omega_X + \ddbar v_{j})^n = \eta_j \omega_X^n,\quad \sup_X v_j = 0. $$
By Theorem \ref{thm:main1} (or \cite{K}), it holds that $\| v_j\|_{L^\infty}\le C_0$ for a uniform $C_0=C_0(n, \omega_X,\epsilon)$.

By integration by parts we have% (c.f. Corollary 3.3 \cite{GZ})
\bea\label{eqn:3.01}
& &\int_X (- u) (\omega_X + \ddbar v_j)^k \wedge \omega_X^{n-k}\\
& \nonumber =& \int_X (-u) \omega_X\wedge \omega_{v_j}^{k-1} \wedge\omega_X^{n-k} + (-v_j) (\omega_u - \omega_X)\wedge \omega_{v_j}^{k-1} \wedge\omega_X^{n-k}\\
& \nonumber \le & \int_X (-u) \omega_X\wedge \omega_{v_j}^{k-1} \wedge\omega_X^{n-k} + C_0\int_X \omega_u \wedge \omega_{v_j}^{k-1} \wedge\omega_X^{n-k}\\
& \nonumber = & \int_X (-u) \omega_{v_j}^{k-1} \wedge\omega_X^{n-k+1} + C_0.%\le C(n) \int_X (-\varphi_t) \omega_X^n + C(n, k, C_0)   \le C,
\eea
Applying (\ref{eqn:3.01}) inductively we get 
\begin{equation}\label{eqn:3.1}
\int_X (- u) (\omega_X + \ddbar v_j)^k \wedge \omega_X^{n-k} \le \int_{X} (-u) \omega_X^n + k C_0\le C,
\end{equation}
for some uniform constant $C = C(n,k,\omega_X,\epsilon)>0$. On the other hand, by Newton-Maclaurin inequality that 
$\big( \frac{\omega_{v_j}^k\wedge \omega_X^{n-k}}{\omega_X^n}\big)^{1/k} \ge c(n,k)\big(\frac{\omega_{v_j}^n}{\omega_X^n}\big)^{1/n}$
we derive from (\eqref{eqn:3.1}) that
$$\int_X (-u) (\eta_j)^{k/n} \omega_X^n \le C(n,k, \omega_X,\epsilon).$$
Letting $j\to\infty $ and applying Fatou's lemma we get
$$\int_K (-u) V_K^{(2\epsilon - 1) k/n} \omega_X^n \le C(n, k,\omega_X,\epsilon)$$ from which we obtain that 
$$V_K = {\mathrm {Vol}} (\{u\le  - s\}) \le C(n, k,\omega_X,\epsilon) s^{-\frac{n}{(2\epsilon - 1) k+ n}}.$$
For any $p< \frac{n}{n-k}$, we have
$$\int_X (-u)^p \omega_X^n \le 1+ p C(n,k, \omega_X,\epsilon) \int_1^\infty s^{p-1 -\frac{n}{(2\epsilon - 1) k+ n}}ds\le C(n, k,\omega_X, p) $$
if $\epsilon = \epsilon(p)>0$ is chosen small enough so that the integral above is integrable. The proof of Lemma \ref{lemma 6} is complete.

\bigskip
Returning to the applications of Theorem \ref{thm:main}, we observe that the condition that $p<{n\over n-k}$ is equivalent to the dual exponent $q$ satisfying $q>{n\over k}$. Thus Theorem \ref{thm:main} combined with Lemmas \ref{eqn:entropy bound} and \ref{lemma new 8} imply at once:

\begin{theorem}
\label{thm:fully}
Consider the family of fully non-linear equations (\ref{eqn:main}) with respect to the degenerating background metrics $\o_t$, $t\in (0,1]$.
Assume that we have solutions $\varphi_t\in C^2$, normalized by $\sup_X\varphi_t=0$. Assume that $\lambda[h_{t,\varphi_t}]\in \Gamma_k$, for some fixed $k$, $1\leq k\leq n$. Fix $q>n/k$.
Then $\|\varphi_t\|_{L^\infty}$ is uniformly bounded by a constant $C$ depending only on 
$n, k, q$ $\omega_X, \chi$  $\| e^{n F_t}\|_{L^q(\omega_X^n)}$ and $\frac{c_t^n}{V_t}$.
\end{theorem}

\medskip

We illustrate this theorem by specializing now to the case of Hessian equations, where $f(\lambda)=\sigma_k(\lambda)^{1\over k}$ for some $1\leq k\leq n$. The more familiar form of this equation is
\begin{equation}\label{eqn:sk section}
(\omega_t + \ddbar \varphi_t)^k\wedge \omega_X^{n-k} = c_t^k e^{k F_t} \omega_X^n,\quad \sup_X \varphi_t = 0,
\end{equation}
and the condition $\lambda[\omega_{t,\varphi_t}]\in \Gamma_k$ is part of the equation\footnote{We remark that usually in the equation (\ref{eqn:sk section}), one normalizes the function $F_t$ such that $\int_X e^{kF_t} \omega_X^n = V$. However, our normalization is that $\int_X e^{n F_t} \omega_X^n = V$.}. Thus Theorem \ref{thm:fully} applies and, assuming uniform bounds for $\|e^{nF}\|_{L^q}$ for some $q>n/k$, it reduces the uniform estimates for $\varphi_t$ to a uniform estimate for the relative volumes $c_t^n/V_t$. An important geometric case when the relative volumes can be controlled is when the classe $\chi$ is big, in the sense that its volume $[\chi^n]=\int_X\chi^n$ is strictly positive. In this case, we obtain

\begin{theorem}
\label{thm:big}
Fix $1\leq k< n$, and consider the family (\ref{eqn:sk section}) of Hessian equations with respect to the degenerating background metrics. Assume that $\chi$ is big. Then for any $q>{n\over k}$, $\|\varphi\|_{L^\infty}$ is bounded uniformly by a constant $C$ depending only on $n, k, q$, $\omega_X, \chi$ and an upper bound for $\| e^{n F_t}\|_{L^q(\omega_X^n)}$.
\end{theorem}

\noindent
{\it Proof of Theorem \ref{thm:big}}. In view of Theorem \ref{thm:fully}, it suffices to show that $c_t^nV_t^{-1}$ is uniformly bounded. Since $V_t\geq [\chi^n]$ for any $t$, this reduces to showing that $c_t$ are themselves uniformly bounded.

The factors $c_t$ are determined by 
\begin{equation}\label{eqn:ct}c_t^k \int_X e^{k F_t} \omega_X^n  = \int_X \omega_t ^k\wedge \omega_X^{n-k}  = O(1).\end{equation} 
To estimate $c_t$, we still need a uniform lower bound of $\int_X e^{kF_t}\omega_X^n$.  We use H\"older's inequality as before.
Thus we write
   \begin{equation}\label{eqn:4.6}
 \int_X e^{k F_t}\omega_X^n\ge \int_{\{F_t\le 0\}} e^{nF_t}\omega_X^n + \int_{\{F_t>0\}} \omega_X^n.
   \end{equation}
Recall that by our normalization on $F_t$ 
   \begin{equation}\label{eqn:4.7}
 V = \int_X e^{n F_t}\omega_X^n =  \int_{\{F_t\le 0\}} e^{nF_t}\omega_X^n + \int_{\{F_t>0\}}e^{n F_t} \omega_X^n.
   \end{equation}
If the first term on the right hand side of (\ref{eqn:4.7}) is greater than $V/2$, then (\eqref{eqn:4.6}) shows that $ \int_X e^{k F_t}\omega_X^n\ge V/2$; otherwise the second term in (\eqref{eqn:4.7}) is greater than $V/2$. Then we have by the H\"older inequality
$$\int_{\{F_t>0\}}e^{n F_t} \omega_X^n \le \Big( \int_X e^{qn F_t} \Big)^{1/q} \Big(\int_{\{F_t>0\}} \omega_X^n  \Big)^{\frac{q - 1}{q}}$$
which yields a uniform lower bound of the second term in (\eqref{eqn:4.6}) depending additionally on the assumed $\| e^{n F_t}\|_{L^q}$. Therefore we conclude from (\eqref{eqn:ct}) that $c_t\le C$ for a uniform constant $C>0$. The proof of Theorem \ref{thm:big} is complete.

\section{Trudinger Inequalities}
\setcounter{equation}{0}

In this section, we illustrate the versatility of our approach by establishing also inequalities of Trudinger type for general non-linear energies. Let $f(\lambda)$ be a fully nonlinear operator satisfying the same structural conditions as in Theorem~\ref{thm:main1}, and define for each $p>0$,
\begin{equation}
E_p(\varphi) = \frac{1}{V}\int_X(-\varphi)^pf^n(\lambda[h_{\varphi}])\omega_X^n
\end{equation}
(this notation is slightly different from the notation $E_t$ used earlier for degenerating background metrics, but there should be no confusion, as the background metric is here fixed, and the index $t$ can be dropped). Recall that we had established before integral estimates for $\varphi$ in terms of the entropy $\Ent_p$. Trudinger estimates are also exponential estimates, but in terms of the energy $E_p$. We have
\begin{theorem}\label{thm: Trudinger-cpt}
	Let $\varphi$ be a plurisubharmonic function such that $\sup \varphi = -1$. Then
	\begin{equation}
	\int_X {\rm exp}\Big\{c(n, p, \gamma)\alpha\Big(\frac{-\varphi}{E_p(\varphi)^{\frac{1}{n+p}}}\Big)^{\frac{n+p}{n}} \Big\}\omega_X^n\leq 2C_{\alpha}
	\end{equation}
	where $\alpha$ and $C_{\alpha}$ are the constants coming from the $\alpha$-invariant estimate of $(X, \omega_X)$. 
\end{theorem}

We note that when specialized to the case when $f(\lambda)$ is the Monge-Amp\`ere operator $f(\lambda) = (\prod_{k=1}^n\lambda_k)^{\frac{1}{n}}$, our theorem recovers an inequality proved in \cite{BB, DGL}. Moreover, our estimate has the major advantage that all constants there depend only on the $\alpha$-invariant of the underlying manifold, hence is uniform over degenerating families with uniform $\alpha$-invariants. 

\medskip
\noindent
{\it Proof of Theorem~\ref{thm: Trudinger-cpt}:}
Let us solve the following auxiliary Monge-Amp\`ere equation with $\sup \psi = 0$, which is solvable due to Yau's theorem \cite{Y}. 
\begin{equation}
(\omega_X+i\p\bar\p\psi)^n = \frac{(-\varphi)^pe^{nF}}{E_p(\varphi)}\omega_X^n
\end{equation}
then by the same maximum argument as in Lemma~\ref{lemma trudinger}, we obtain the inequality 
\begin{equation}\label{eq: almost-Trudinger}
{c(n, p, \gamma)  \Big(\frac{-\varphi}{E_p(\varphi)^{1/(n+p)}}\Big)^{\frac{n+p}{n}}}\leq -\psi+C(n, p, \gamma)E_p(\varphi)^{\frac{1}{p}}
\end{equation}
Now we pick $\kappa = 2^{\frac{n}{n+p}}C^{\frac{n}{n+p}}c^{-\frac{n}{n+p}}$, where $C$ and $c$ are constants in estimate above, which only depends on $n, p$ and $\gamma$. Now set $U_{\kappa} = \{-\varphi\leq \kappa E_p(\varphi)^{\frac{1}{p}}\}$. Then on $X\setminus U_{\kappa}$, we have by our choice of $\kappa$,
\begin{equation}
\frac{1}{2}c(n, p, \gamma)\Big(\frac{-\varphi}{E_p(\varphi)^{1/(n+p)}}\Big)^{\frac{n+p}{n}}\leq -\psi
\end{equation}
and on $U_{\kappa}$, we have 
\begin{equation}
c(n, p, \gamma)\Big(\frac{-\varphi}{E_p(\varphi)^{1/(n+p)}}\Big)^{\frac{n+p}{n}} \leq -c(n, p, \gamma)\kappa^{\frac{p}{n}}\varphi 
\end{equation}
now multiplying by $\min(c^{-1}\kappa^{-\frac{p}{n}}, 1/2)\alpha$ and integrating, we get
\begin{equation}
\int_X{\rm exp}\Bigg\{c'(n, p, \gamma)\alpha \Big(\frac{-\varphi}{E_p(\varphi)^{1/(n+p)}}\Big)^{\frac{n+p}{n}}\Bigg\}\omega_X^n \leq \int_{U_{\kappa}}e^{-\alpha\psi}\omega_X^n+\int_{X\setminus U_{\kappa}}e^{-\alpha\varphi}\omega_X^n\leq 2C_{\alpha}
\end{equation}
which is the desired result.

\bigskip

%Email address:  phong@math.columbia.edu, tong@math.columbia.edu

\noindent Department of Mathematics \& Computer Science, Rutgers University, Newark, NJ 07102, USA

\noindent bguo@rutgers.edu,

\medskip

\noindent Department of Mathematics, Columbia University, New York, NY 10027 USA

\noindent phong@math.columbia.edu, tong@math.columbia.edu

\end{document}